\documentclass[a4paper]{amsart}
\usepackage{amsfonts}
\theoremstyle{plain}
\begingroup
\newtheorem{theorem}{Theorem}[section]
\newtheorem{lemma}[theorem]{Lemma}
\newtheorem{proposition}[theorem]{Proposition}
\newtheorem{corollary}[theorem]{Corollary}
\newtheorem{remark}[theorem]{Remark}

\theoremstyle{definition}
\theoremstyle{remark}
\numberwithin{equation}{section}

\newcommand{\hs}{{\mathcal H}}
\newcommand{\ks}{{\mathcal K}}
\newcommand{\cs}{{\mathcal C}}
\newcommand{\gs}{{\mathcal G}}
\newcommand{\fs}{{\mathcal F}}
\newcommand{\leb}{{\mathcal L}}

\newcommand{\ms}{{\mathcal M}}

\newcommand{\eub}{{\mathcal E}}
\newcommand{\vub}{{\mathcal V}}

\newcommand{\R}{{\mathbb R}}

\newcommand{\msim}{{\rm M}^{2\times 2}_{\rm sym}}

\newcommand{\Ba}{B_1(0)}
\newcommand{\Bb}{\overline{B}_1(0)}

\newcommand{\Om}{\Omega}
\newcommand{\Omb}{\overline{\Omega}}

\newcommand{\deli}[2]{L^{1,2}(#1 \setminus #2)}

\newcommand{\Ni}{\Gamma}

\newcommand{\nL}[1]{||\nabla #1||_a}
\newcommand{\nA}[1]{||E #1||_A}

\newcommand{\weakst}{\stackrel{\ast}{\rightharpoonup}}
\newcommand{\weakstloc}{{\stackrel{\ast}{\rightharpoonup}}_{loc}}
\newcommand{\weak}{\rightharpoonup}

\newcommand{\esbd}[1]{\partial^* \!#1}
\newcommand{\per}[2]{{\rm P}(#1,#2)}

\newcommand{\fun}[1]{\int_{#1} \varphi(x,\nu) \,d{{\mathcal H}}^1}
\newcommand{\funx}[1]{\int_{#1} \varphi(x,\nu_x) \,d{{\mathcal H}}^1}
\newcommand{\funU}[1]{\int_{#1 \cap U} \varphi(x,\nu) \,d{{\mathcal H}}^1}
\newcommand{\pr}[2]{(\nabla #1,\nabla #2)_a}
\newcommand{\prA}[2]{(E #1,E #2)_A}

\newcommand{\res}{\mathop{\hbox{\vrule height 7pt width .5pt depth 0pt
\vrule height .5pt width 6pt depth 0pt}}\nolimits}

\title
[A GENERALIZATION OF GO\L\c{A}B'S THEOREM]
{A GENERALIZATION OF GO\L\c{A}B'S THEOREM \\
AND \\
APPLICATIONS TO FRACTURE MECHANICS}
\author[A. Giacomini]
{Alessandro Giacomini}
\address[Alessandro Giacomini]{S.I.S.S.A., Via Beirut 2-4, 34014, Trieste,
Italy}
\email[A. Giacomini]{giacomin@sissa.it}
\begin{document}
\vskip .2truecm
\begin{abstract}
\small{We study the lower semicontinuity for functionals
of the form $K \to \fun{K}$ defined on compact sets in $\R^2$
with a finite number of connected components and
finite $\hs^1$ measure and apply the result to the study of quasi-static growth of
brittle fractures in linearly elastic inhomogeneous and anisotropic bodies.}
\end{abstract}
\maketitle

\section{Introduction}
\label{intr}

In 1998, G.A. Francfort and J.-J. Marigo
\cite{FM} proposed a model for the quasi-static growth of
brittle fractures in elastic bodies.
This model is based on Griffith's criterion of crack growth
which takes into account a competition between the bulk energy
given by the deformation and the surface energy given by the
length of the fracture.
Recently, G. Dal Maso and R. Toader \cite{DMT} gave a precise
mathematical formulation of the model in  {\it dimension
two} for linearly elastic homogeneous bodies under anti-planar
shear.
\par
The aim of this paper is to extend
this analysis in {\it dimension two} to
{\it anisotropic linearly elastic inhomogeneous
bodies} subjected to {\it anti-planar} or {\it planar shear}.
{\it Anisotropy} will be considered both in the {\it bulk}
and in the {\it surface energy}.
\par
In order to make the ideas precise, let $\Om \subseteq \R^2$ be open and
bounded and let $\ks^f_m(\Omb)$ denote the family of compact subsets
of $\Omb$ with at most $m$ connected components and
finite $\hs^1$ measure .
Consider an elastic body
of the form $\Om \times \R$ and assume the cracks of the form
$K \times \R$ with $K \in \ks^f_m(\Omb)$. Assume the
displacement $u\,:\,\Om \times \R \to \R^3$ depends only on
$x_1, x_2$. If $u(x_1,x_2)=(0,0,u_3(x_1,x_2))$, we are in the
case of {\it anti-planar shear}, while if
$u(x_1,x_2)=(u_1(x_1,x_2),u_2(x_1,x_2),0)$, we speak of
{\it planar shear}.
In the model case, we consider the {\it bulk energy} (referred to
a finite portion of the cylinder determined by two cross
sections separated by a unit distance) of the form
$$
\int_{\Om \setminus K} \mu |Eu|^2+ \lambda |trEu|^2 \,d\leb^2
$$
where $\mu, \lambda$ are called {\it Lam\'e
coefficients} and $Eu$ is
the symmetric part of the gradient of $u$.
\par
The {\it surface energy}
on the fracture $K$ is given by
$$
\int_{K} \varphi(x,\nu_x)\,d\hs^1(x),
$$
where $\nu_x$ is the
{\it unit normal vector} at $x$ to $K$ and
$\varphi\,:\, \Omb \times \R^2 \to [0,\infty[$
is a continuous function, positively $1$-homogeneous, even and
convex in the second variable.
\par
Given $\partial_D \Om \subseteq \partial \Om$ open in the
relative topology and with a finite number of connected
components, we prescribe a displacement $g$ on $\partial_D \Om$.
The displacement $u_{g,K}$ of the elastic body relative to $g$
and the crack $K$ is obtained minimizing
$\int_{\Om \setminus K} \mu |Eu|^2+ \lambda |trEu|^2 \,d\leb^2$
under the condition $u=g$ on $\partial_D \Om \setminus K$.
The condition $u=g$ on $\partial_D \Om \setminus K$
takes into account the fact that the displacement is not
transmitted through a fractured region.
The {\it total elastic energy} is
given by
$$
\eub(g,K):=
\int_{\Om \setminus K}
\mu |Eu_{g,K}|^2+ \lambda |trEu_{g,K}|^2\,d\leb^2
+ \fun{K}.
$$
\par
Suppose an initial crack $K_0$ and boundary displacements $g(t)$,
$t\in[0,1]$, $g(0)=0$, are given. By a {\it quasi-static growth} of the
fracture, we
mean an {\it increasing} map $t \to K(t)$ from $[0,1]$ to $\ks^f_m(\Omb)$
with $K(0)=K_0$ and such that $K(t)$ minimizes $\eub(g(t),K)$
among $K$'s such that $\cup_{s<t}K(s) \subseteq K$.
The constraints given by the previous cracks
indicate the {\it irreversibility} of the growth and the {\it absence of
healing phenomena}.
We require also the {\it stationarity condition}
$\frac{d}{ds} \eub(g(t),K(s))|_{s=t}=0$ and the {\it absolute
continuity} of the total
elastic energy $t \to \eub(g(t),K(t))$ even if, as noted in
\cite{FM}, $t \to \fun{K(t)}$ could be discontinuous.
\par
The quasi-static growth
satisfying the stationarity condition and the absolute continuity
of the total energy
is obtained through
a {\it time discretization method}. Given $\delta>0$, we divide
$[0,1]$ in $N_\delta$ intervals $[t^\delta_i,t^\delta_{i+1}]$
and we indicate by $K^\delta_i$ the solution of
$$
\min \{ \eub(g(t^\delta_i),K)\,:\,
K^\delta_{i-1} \subseteq K \},
$$
where we consider $K^\delta_{-1}=K_0$. We make the interpolation
$K^\delta(t)=K^\delta(t^\delta_i)$ if
$t^\delta_i \le t < t^\delta_{i+1}$; letting $\delta \to 0$ along
a suitable sequence, it turns out that
$K^\delta(t) \to K(t)$ in the Hausdorff metric determining the
quasi-static growth. Moreover
$Eu^\delta(t) \to Eu(t)$ strongly, so that the elastic bulk
energies of the approximating fractures converge to the bulk
energy of the solution.
\par
Moreover we will prove
that
$\fun{K^\delta(t)}$ converges to the surface energy of the
solution. We conclude that the time discretization
procedure gives an approximation both of the bulk and surface
energy of the solution. We remark that this fact
is new also in the case $\varphi \equiv 1$, that is
when the surface energy depends only on the length of the fracture.
\par
In order to deal with an {\it anisotropic} and {\it inhomogeneous}
surface energy, the main step is to prove a lower semicontinuity
theorem for the functional $\fs(K):= \funx{K}$ on $\ks^f_m(\Omb)$ with
respect to Hausdorff convergence.
This functional is well defined: in fact,
even if $K$ is
not in general the union of $m$ regular curves, it
turns out that it is possible
to define at $\hs^1$-a.e. point $x \in K$ an {\it approximate unit
normal vector} $\nu_x$ completely determined up to the sign. In
the case $K$ is regular, $\nu_x$ coincides with the usual normal
vector.
Note that for $\varphi \equiv 1$,
the semicontinuity result reduces to Go\l\c{a}b's theorem on
the lower semicontinuity of
$\hs^1$ measure under Hausdorff convergence.
\par
The paper is organized as follows.
After some preliminaries, we prove the lower semicontinuity result
in Section \ref{sem}. In sections \ref{antplan} and \ref{plan},
we deal with the study of quasi-static growth of brittle
fractures in the anti-planar and planar
cases. Using shape continuity results proved in \cite{Buc} and
\cite{Ch}, we can treat inhomogeneous bulk energies:
we consider quadratic forms of $Eu$ equivalent to the standard
$\int_{\Om \setminus K} |Eu|^2\,d\leb^2$.
This cannot be done directly using the techniques
of \cite{DMT} where the strong convergence of the gradient of
deformation is obtained through a duality argument which relies
on the particular form $\int_{\Om \setminus K} |\nabla u|^2\,d\leb^2$
of the bulk energy.

\section{ NOTATIONS AND PRELIMINARIES}
\label{prel}

In what follows, $\Om \subseteq \R^2$ is a bounded open set
with Lipschitz boundary, $\partial_D \Om$ is a subset of $\partial \Om$
open in the relative topology and with a finite number of connected components.
\bigskip
\par\noindent
{\it Sets with finite perimeter.}
We indicate the perimeter of $E$ in $\Om$ by $\per{E}{\Om}$. Let $E$ be a set of finite perimeter
in $\Om$; the reduced boundary $\esbd{E}$ and the approximate inner normal $\nu$ at points of
$\esbd{E}$ are defined such that the following identity holds:
$$
\forall g \in C_c(\Om,\R^2)\,\,
-\int_E {\rm div} g \,d\leb^2= \int_{\esbd{E}} g \cdot \nu d\,\hs^1.
$$
Set $\mu_E= \nu \hs^1 \res \esbd{E}$. For all $x \in \esbd{E}$, indicated
the map $\xi \to \frac{1}{\lambda} (\xi-x)$ by $D_\lambda$, the following
{\it blow up property} holds: for $\lambda \to 0^+$
$$
\mu_{D_\lambda(E)} \weakstloc \hs^1 \res T_\nu,
$$
locally weakly star in the sense of measures, where $T_\nu$ is the subspace
orthogonal to $\nu$.
\par
We say that a sequence $(E_h)$ of subset of $\Om$ converges in $L^1_{\rm loc}(\Om)$ to $E$,
if the corresponding characteristic functions $\chi_{E_h}$ converge in $L^1_{\rm loc}(\Om)$ to
$\chi_E$. If there exists $C \ge 0$ such that $\per{E_h}{\Om} \le C$ for all $h$
and $E_h \to E$ in $L^1_{\rm loc}(\Om)$, then $E$ has finite perimeter in $\Om$ and
$\mu_{E_h} \weakst \mu_E$ in the weak star topology of $\ms_b(\Om,\R^2)$.
For further details on sets of finite perimeter, the reader is referred to \cite{AFP}.
\vskip 20pt
\noindent
{\it Hausdorff metric on compact sets.}
We indicate the set of all compact subsets of $\Omb$ by $\ks(\Omb)$,
the set of elements of $\ks(\Omb)$
with finite $\hs^1$ measure by $\ks^f(\Omb)$ and, given $\lambda \ge 0$, the compact
sets $K$ with $\hs^1(K) \le \lambda$ by $\ks^\lambda(\Omb)$.
$\ks(\Omb)$ can be endowed by the
Hausdorff metric $d_H$ defined by
$$
d_H(K_1,K_2) := \max \left\{ \sup_{x \in K_1} {\rm dist}(x,K_2), \sup_{y \in
K_2} {\rm dist}(y,K_1)\right\}
$$
with the conventions ${\rm dist}(x, \emptyset)= {\rm diam}(\Om)$ and $\sup
\emptyset=0$, so that $d_H(\emptyset, K)=0$ if $K=\emptyset$ and
$d_H(\emptyset,K)={\rm diam}(\Om)$ if $K \not=\emptyset$. It turns out that
$\ks(\Omb)$ endowed with the Hausdorff metric is a compact space
(see e.g. \cite{Ro}).
Let $\ks_m(\Omb)$ be the subset of $\ks(\Omb)$ of those compact sets
which have less than $m$ connected components.
Since Hausdorff convergence preserves connectedness, $\ks_m(\Omb)$
are closed subsets of $\ks(\Omb)$ for all $m$. Let
$\ks_m^f(\Omb):=\ks_m(\Omb) \cap \ks^f(\Omb)$ and given $\lambda \ge 0$,
$\ks_m^\lambda(\Omb):=\ks_m(\Omb) \cap \ks^\lambda(\Omb)$.
\par
Hausdorff measure $\hs^1$ is not lower semicontinuous in $\ks(\Omb)$ with respect to
Hausdorff metric. However it is lower semicontinuous if restricted to
$\ks_m(\Omb)$: for the case $m=1$, this result is known as Go\l\c{a}b's theorem
(see e.g. \cite{MS}). The general case can be found in \cite{DMT}.

\begin{theorem}
\label{golab}
Let $(K_n)$ be a sequence in $\ks_m(\Omb)$ which
converges to $K$ in the Hausdorff metric. Then $K \in \ks_m(\Omb)$ and for every open
subset $U \subseteq \R^2$
$$
\hs^1(K \cap U) \le \liminf_n \hs^1(K_n \cap U).
$$
\end{theorem}

\vskip20pt
\noindent
{\it Structure of compact connected sets with finite $\hs^1$ measure.}
\;It can be proved (see e.g. \cite{Fal}) that
if $K \in \ks^f_1(\Omb)$, for a.e. $x \in K$ there exists an approximate unit normal
vector $\nu_x$ which is characterized by
\begin{equation}
\label{blowup}
\mu_{D_\lambda(E)} \weakstloc \hs^1 \res T_{\nu_x}
\quad {\rm for}\; \lambda \to 0^+
\end{equation}
locally weakly star in the sense of measures, where $T_{\nu_x}$ is the
subspace of $\R^2$ orthogonal to $\nu_x$.
Moreover the map ${x \to \nu_x}$ is Borel measurable, so that for every
continuous function
$\varphi:\Omb \times \R^2 \to [0,\infty[$
even in the second variable the integral
$$
\fun{K}
$$
is well defined. Clearly the functional is well defined also for $K \in \ks^f_m(\Omb)$ with
$m \ge 1$.
\par
In section \ref{sem} we will be concerned in the problem of the lower semicontinuity of
the function $K \to \fun{K}$ under the Hausdorff convergence.
\par
We will use the fact that a connected set $C$ with finite $\hs^1$ measure is arcwise connected
and moreover $\hs^1(\overline{C})=\hs^1(C)$: see e.g. \cite{DMT}.
\vskip 20pt
\noindent
{\it Reshetnyak's theorems on measures.} The following theorem gives a
lower semicontinuity result for functionals defined on measures; for a proof, the
reader is referred to \cite{AFP}. If $\mu$ is a measure, let $|\mu|$ be its
total variation and let $\frac{d\mu}{d|\mu|}$ be the Radon-Nicodym derivative of $\mu$
with respect to $|\mu|$.
\begin{theorem}
\label{Resh}
Let $\Om$ be an open subset of $\R^n$ and $\mu,\mu_k$ be $\R^m$-valued finite
Radon measures
in $\Om$; if $\mu_h \to \mu$ weakly star in $\ms_b(\Om,\R^m)$ then
$$
\int_\Om f \left( x, \frac{d\mu}{d|\mu|}(x) \right) d|\mu|(x) \le
\liminf_{h \to \infty}
\int_\Om f \left( x, \frac{d\mu_h}{d|\mu_h|}(x) \right) d|\mu_h|(x)
$$
for every lower semicontinuous function $f:\,\Om \times \R^m \to [0,+\infty]$, positively
$1$-homogeneous and convex in the second variable.
\end{theorem}
We say that $\mu_n$ converges strictly to $\mu$ in $\ms_b(\Om,\R^m)$ if $\mu_n \to \mu$ weakly star
and $|\mu_n|(\Om) \to |\mu|(\Om)$.
The following theorem gives a continuity result for functional defined on
measures: for a proof
see \cite{AFP}.

\begin{theorem}
\label{Resh2}
Let $\Om$ be an open subset of $\R^n$ and $\mu,\mu_k$ be $\R^m$-valued finite
Radon measures
in $\Om$; if $\mu_h \to \mu$ strictly in $\ms_b(\Om,\R^m)$ then
$$
\lim_{h \to \infty}
\int_\Om f \left( x, \frac{d\mu_h}{d|\mu_h|}(x) \right) d|\mu_h|(x)=
\int_\Om f \left( x, \frac{d\mu}{d|\mu|}(x) \right) d|\mu|(x)
$$
for every continuous and bounded function $f:\,\Om \times S^{m-1} \to \R$.
\end{theorem}

\vskip20pt
\noindent
{\it Deny-Lions spaces.}
If $A$ is an open subset of $\R^2$, the Deny-Lions space $L^{1,2}(A)$ is defined as
\begin{equation}
\label{DL}
L^{1,2}(A):=
\left\{ u \in W^{1,2}_{\rm loc}(A)\,:\, \nabla u \in L^2(A, \R^2)
\right\}.
\end{equation}
In the case in which $A$ is regular $L^{1,2}(A)$ coincides with the usual Sobolev space while
if it is irregular, it can be strictly larger.
In what follows, given $K \subseteq \Omb$ compact and
$u \in L^{1,2}(\Om \setminus K)$, we extend $\nabla u$ to $0$ on $K$,
so that $\nabla u \in L^2(\Om, \R^2)$ although $\nabla u$ is the distributional derivative of
$u$ only on $\Om \setminus K$.
The following theorem proved in \cite{Buc} will be used in Section \ref{antplan}.

\begin{theorem}
\label{p2}
Let $m \ge 1$, $K_n$ a sequence in $\ks_m(\Omb)$ which
converges to $K$ in the Hausdorff metric and such that
$\leb^2(\Om \setminus K_n) \to \leb^2(\Om \setminus K)$.
Then for every $u \in L^{1,2}(\Om \setminus K)$,
there exists $u_n \in  L^{1,2}(\Om \setminus K_n)$ such that $\nabla u_n \to \nabla u$ strongly
in $L^2(\Om,\R^2)$.
\end{theorem}

Consider now for $A$ open subset of $\R^2$
\begin{equation}
\label{LD}
LD^{1,2}(A):=
\left\{ u \in W^{1,2}_{\rm loc}(A;\R^2)\,:\, E(u) \in L^2(A,\msim)
\right\},
\end{equation}
where
$Eu:=\frac{1}{2}(\nabla u+ (\nabla u)^t)$
is the symmetric part of the gradient of $u$
and $\msim$ is the space of $2 \times 2$ symmetric matrices endowed with the standard scalar product
$B_1{:}B_2:= tr(B_1^tB_2)$ and the corresponding norm $|B|:=(B{:}B)^{\frac{1}{2}}$.
\par
In what follows, given $K \subseteq \Omb$ compact and
$u \in LD^{1,2}(\Om \setminus K)$,
we extend $Eu$ to $0$ on $K$ although it coincides with the symmetric part of the
distributional gradient of $u$ only on $\Om \setminus K$.
The following result, which can be obtained combining the density result proved in
\cite{Ch} and Theorem \ref{p2}, will be used in Section \ref{plan}.

\begin{theorem}
\label{Cham}
Let $m \ge 1$, $K_n$ a sequence in $\ks_m(\Omb)$ which
converges to $K$ in the Hausdorff metric and such that
$\leb^2(\Om \setminus K_n) \to \leb^2(\Om \setminus K)$.
Then for every $u \in LD^{1,2}(\Om \setminus K)$, there exists a sequence
$u_n \in LD^{1,2}(\Om \setminus K_n)$ such that $Eu_n \to Eu$ strongly in $L^2(\Om; \msim)$.
\end{theorem}

\vskip20pt
\noindent
{\it Absolutely continuous function.}
Given a Hilbert space $X$, we indicate by $AC([0,1],X)$ the space of absolutely continuous
function from $[0,1]$ to $X$: for the main properties of this space, the reader is referred
to \cite{Br1}. Given $g \in AC([0,1],X)$, the time derivative, which exists a.e. in $[0,1]$,
is denoted by $\dot{g}$.

\section{THE MAIN RESULTS}
\label{main}

Let $\varphi\,:\,\Omb \times \R^2 \to [0,+\infty[$ a continuous function,
positively 1-homogeneous, even and convex in the second variable such
that for $c_1,c_2>0$
\begin{equation}
\label{varphiprop}
\forall (x,\nu) \in \Omb \times \R^2 \,:\, c_1|\nu| \le \varphi(x,\nu) \le c_2|\nu|.
\end{equation}

The main result of the paper is the following lower semicontinuity theorem.

\begin{theorem}
\label{main0}
The functional
$$
\begin{array}{r}
\fs\,:\,\ks_m^f(\Omb) \\  \\ K
\end{array}
\begin{array}{c}
\longrightarrow \\  \\ \mapsto
\end{array}
\begin{array}{l}
[0,{\infty}{[} \\ \\ \fun{K}
\end{array}
$$
is lower semicontinuous if $\ks_m^f(\Omb)$ is endowed with the Hausdorff metric.
\end{theorem}

The previous theorem will be used to deal with the problem of evolution of
brittle fractures in linearly elastic bodies.
\par
Let $a \in L^\infty(\Om,\msim)$ such that for $\alpha_1, \alpha_2 >0$
\begin{equation}
\label{unifellipt}
\forall x \in \Om \,,\forall \xi \in \R^2\,:\,
\alpha_1 |\xi|^2 \le a(x)\xi \cdot \xi \le \alpha_2 |\xi|^2.
\end{equation}
Let $(\cdot,\cdot)_a$ denote the associated scalar product on $L^2(\Om, \R^2)$
defined as
$$
(v,w)_a= \int_\Om \sum_{i,j=1}^2 a(x) v_i(x) w_j(x) \,d\leb^2(x)
$$
and let $||\cdot||_a$ be the relative norm.
\par
For every $g \in H^1(\Om)$ and $K \in \ks_m^f(\Omb)$, we set
\begin{equation}
\label{prantplan}
\eub(g,K):= \min_{v \in \Ni(g,K)}
\left\{ \int_{\Om \setminus K} a(x) \nabla v \cdot \nabla v \,d\leb^2 + \fun{K} \right\},
\end{equation}
where
\begin{equation}
\label{gamma}
\Ni(g,K) :=\left\{ u \in \deli{\Om}{K}\,,\, u=g\,{\rm on}\, \partial_D \Om \setminus K \right\}.
\end{equation}
\par
The following theorem states the existence of a quasi-static evolution for brittle fractures in linear elastic bodies
under anti-planar displacement: note that
both the bulk and the surface energy depend in a possibly inhomogeneous way
on the anisotropy of the body.

\begin{theorem}
\label{main1}
Let $m \ge 1$, $g \in AC([0,1],H^1(\Om))$, $K_0 \in \ks_m^f(\Omb)$.
Then there exists a function $K:[0,1] \to \ks_m^f(\Omb)$ such that, letting
$u(t)$ be a solution of the minimum problem
(\ref{prantplan}) which defines $\eub(g(t),K(t))$ for all $t \in [0,1]$,
\begin{itemize}
\item[]
\item[(a)] $\hfill \displaystyle K_0 \subseteq K(s) \subseteq K(t) \; {\it for} \; 0\le s \le t \le 1; \hfill$
\item[]
\item[(b)] $\hfill \displaystyle \eub(g(0),K(0)) \le \eub(g(0),K) \quad \forall K \in \ks_m^f(\Omb),\, K_0 \subseteq K;\hfill$
\item[]
\item[(c)] $\hfill \displaystyle \forall t \in ]0,1]\,:\,\eub(g(t),K(t)) \le \eub(g(t),K) \quad
\forall K \in \ks_m^f(\Omb),\, \cup_{s<t} K(s) \subseteq K;\hfill$
\item[]
\item[(d)] $\hfill \displaystyle t \mapsto \eub(g(t),K(t))\, {\it is\,absolutely\,continuous\,on} \,[0,1];\hfill$
\item[]
\item[(e)] $\hfill \displaystyle \frac{d}{dt} \eub(g(t),K(t))=2\pr{u(t)}{\dot{g}(t)}
\quad {\it for\,a.e.}\; t\in [0,1],\hfill$
\item[]
\item[(f)] $\hfill \displaystyle \frac{d}{ds} \eub(g(t),K(s)) \vert_{s=t}=0 \quad {\it for \;a.e.} \;t \in [0,1]. \hfill$
\end{itemize}
\end{theorem}

\bigskip
Let $\leb(\msim)$ be the space of automorphism of $\msim$
and let
$A \in L^\infty (\Om,\leb(\msim))$ such that there exist
$\alpha_1, \alpha_2 >0$ with
$$
\forall x \in \Om \,:\, \alpha_1 |M|^2 \le A(x)M{:}M \le \alpha_2 |M|^2.
$$
Let us pose
$\prA{u}{v}:= \int_{\Om \setminus K} A(x)Eu{:}Ev\,d\leb^2$ and
$\nA{u}:=\prA{u}{u}^{\frac{1}{2}}$.
\par
For every $g \in H^1(\Om;\R^2)$ and $K \in \ks_m^f(\Omb)$, set
\begin{equation}
\label{prplan}
\gs(g,K)= \min_{v \in \vub(g,K)}
\left\{ \int_{\Om \setminus K} A(x)Eu{:}Eu\,d\leb^2+ \fun{K} \right\},
\end{equation}
where
\begin{equation}
\label{ni}
\vub(g,K) =\left\{ u \in
LD^{1,2}(\Om \setminus K)\,,\, u=g\,{\rm on}\, \partial_D \Om \setminus K \right\}.
\end{equation}
The following theorem states the existence of a quasi-static evolution for
brittle fractures in inhomogeneous anisotropic linearly elastic bodies
under planar displacement.

\begin{theorem}
\label{main2}
Let $m \ge 1$, $g \in AC([0,1],H^1(\Om;\R^2))$, $K_0 \in \ks_m^f(\Omb)$.
Then there exists a function $K:[0,1] \to \ks_m^f(\Omb)$ such that, letting
$u(t)$ be a solution of the minimum problem
(\ref{prplan}) which defines $\gs(g(t),K(t))$ for all $t \in [0,1]$,
\begin{itemize}
\item[]
\item[(a)] $\hfill \displaystyle K_0 \subseteq K(s) \subseteq K(t) \; {\it for} \; 0\le s \le t \le 1; \hfill$
\item[]
\item[(b)] $\hfill \displaystyle \gs(g(0),K(0)) \le \gs(g(0),K) \quad \forall K \in \ks_m^f(\Omb),\, K_0 \subseteq K;\hfill$
\item[]
\item[(c)] $\hfill \displaystyle \forall t \in ]0,1]\,:\,\gs(g(t),K(t)) \le \gs(g(t),K) \quad
\forall K \in \ks_m^f(\Omb),\, \cup_{s<t} K(s) \subseteq K;\hfill$
\item[]
\item[(d)] $\hfill \displaystyle t \mapsto \gs(g(t),K(t))\, {\it is\,absolutely\,continuous\,on} \,[0,1];\hfill$
\item[]
\item[(e)] $\hfill \displaystyle \frac{d}{dt} \gs(g(t),K(t))=2\prA{u(t)}{\dot{g}(t)}
\quad {\it for\,a.e.}\; t\in [0,1],\hfill$
\item[]
\item[(f)] $\hfill \displaystyle \frac{d}{ds} \gs(g(t),K(s)) \vert_{s=t}=0 \quad {\it for \;a.e.} \;t \in [0,1]. \hfill$
\end{itemize}
\end{theorem}

\begin{remark}
\label{mainrem}
{\rm
It turns out that for every function $K:[0,1] \to \ks^f_m(\Omb)$ which satisfies (a)-(d) of
Theorem \ref{main1}, then conditions (e) and (f) are equivalent.
Similarly, for every function $K:[0,1] \to \ks^f_m(\Omb)$ which satisfies (a)-(d) of Theorem \ref{main2}, conditions (e)
and (f) are equivalent.
}
\end{remark}

We will prove theorem \ref{main0} in section \ref{sem} using a comparison of measures which involves a blow-up
technique; theorems \ref{main1} and \ref{main2} will be proved in sections \ref{antplan} and \ref{plan}
respectively: a discretization in time procedure will be employed and,
in the particular case in which $g(0)=0$, we prove that
this method gives an approximation of the total energy of the solution.

\section{A GENERALIZATION OF GO\L\c{A}B THEOREM}
\label{sem}
Throughout this section, let $\varphi: \Omb \times \R^2 \to [0, \infty[$ be a
continuous function, positively $1$-homogeneous, even and convex in
the second variable satisfying
\begin{equation}
\label{eq1}
\forall \nu \in \R^2\,:\,c_1 |\nu| \le \varphi(x,\nu) \le c_2 |\nu|
\end{equation}
for some $c_1,c_2 >0$.
\par
Let $\cs$ be the subset of $L^1(\Om)$ composed by characteristic functions of sets
with finite perimeter in $\Om$.

\begin{theorem}
\label{sem.1}
Consider the functional $\gs: \cs \to [0,\infty[$ defined by
$$
\gs(E)= \fun{\esbd{E}}
$$
where $\nu$ denotes the inner normal of $E$. Then $\gs$ is lower semicontinuous
with respect to the $L^1$ topology.
\end{theorem}

\begin{proof}
Let $(E_h)$ be a sequence of sets with finite perimeter in $\Om$ with $E_h \to E$ in
$L^1(\Om)$: it is sufficient to consider the case $\per{E_h}{\Om} \le C$ for some $C \ge 0$
independent of $h$.
As noted in Section \ref{prel}, $\mu_{E_h} \weakst \mu_E$ in the weak star topology of
$\ms_b(\Om,\R^2)$.
Since the inner normal to $E_h$ (resp. $E$) is given by $\displaystyle \frac{d\mu_{E_h}}{d\hs^1}$
(resp. by $\displaystyle \frac{d\mu_{E}}{d\hs^1}$), we can use Reshetnyak lower semicontinuity theorem (see
Section \ref{prel}) to get the conclusion.
\end{proof}

\begin{theorem}
\label{sem.2}
Let $U$ be an open subset of $\R^2$. The functional
$$
\begin{array}{r}
\fs\,:\,\ks_m^f(\Omb) \\  \\ K
\end{array}
\begin{array}{c}
\longrightarrow \\  \\ \mapsto
\end{array}
\begin{array}{l}
[0,{\infty}{[} \\ \\ \funU{K}
\end{array}
$$
is lower semicontinuous if $\ks_m^f(\Omb)$ is endowed with the Hausdorff metric.
\end{theorem}

\begin{proof}
We consider preliminarily the case $m=1$.
\par
Let $K_n,K \in \ks^f_1(\Omb)$ with $K_n \to K$ in the Hausdorff metric: our aim is to verify that
$$
\funU{K} \le \liminf_n \funU{K_n}.
$$
Without loss of generality we may consider sequences $(K_n)$ such that
$$
\sup_n \funU{K_n} < +\infty.
$$
Let us consider the positive measures $\mu_n, \mu$ in $\ms_b(U)$
$$
\mu_n(B)= \fun{K_n \cap B},
$$
$$
\mu(B)= \fun{K \cap B}.
$$
By (\ref{eq1}), $(\mu_n)$ is bounded in $\ms_b(U)$ and so up to a
subsequence it converges in the weak-star topology of $\ms_b(U)$ to a measure $\mu_0$ whose
support is contained in $K \cap U$. By weak convergence we have
$$
\mu_0(U) \le \liminf_n \mu_n(U),
$$
and so it is sufficient to prove that
\begin{equation}
\label{eqcfr}
\mu(U) \le \mu_0(U).
\end{equation}
We prove instead that $\mu \le \mu_0$ using a density
argument which requires a blow-up technique: we obtain (\ref{eqcfr}) as a consequence.
\par
Firstly consider $K_n \in \ks_1^f(\Bb)$, $\hs^1(K_n) \le C$
for some $C \ge 0$ and $K_n \to K$ in the Hausdorff metric where $K$
is the diameter connecting the points $e_1:=(1,0)$ and $-e_1$.
Note that for every strip $S_\eta=\{x\in \R^2\,:\, -\eta \le x_2 \le \eta\}$ with $\eta>0$,
$K_n \subseteq S_\eta$ and $K_n \cap \partial S_\eta = \emptyset$ for $n$ large enough.
Given $\varepsilon>0$, let $V^\varepsilon:=\{x\in \R^2\,:\, -1+\varepsilon \le x_1 \le 1-\varepsilon\}$,
$\partial^{\pm}V^\varepsilon$
the connected components of $\partial V^\varepsilon$ containing the points $(1-\varepsilon)e_1$ and
$-(1-\varepsilon)e_1$ respectively.
For $n$ large enough, since $K_n$ is connected, there exist points
$x^{\pm}_n \in \partial^{\pm}V^{\varepsilon} \cap K_n$ such that
$x^{\pm}_n \to \pm (1-\varepsilon)e_1$. Let $L_n$ be the union of the segments connecting $x_n^-$ to
$-(1-\varepsilon)e_1$ and $-(1-\varepsilon)e_1$ to $-e_1$, $x_n^+$ to
$(1-\varepsilon)e_1$ and $(1-\varepsilon)e_1$ to $e_1$. Note that $H_n:=K_n \cup L_n$ is connected and that
$$
\hs^1(L_n) \le 3\varepsilon
$$
for $n$ large enough.
\par
Let $E_n$ be the connected component of $\Ba \setminus H_n$
containing $\frac{1}{2}e_2$, where $e_2:=(0,1)$. As $\pm e_1 \in H_n$ and $H_n$ converges to $K$
in the Hausdorff metric, it is easy to see that $E_n$ converges in $L^1$ to
$B_1^+(0):=\{x \in B_1(0)\,:\, x_2>0\}$.
$E_n$ has finite perimeter because $\partial E_n \subseteq H_n$ and
these sets have finite $\hs^1$ measure (see Proposition
3.62 of \cite{AFP}). By Theorem \ref{sem.1} we have
\begin{eqnarray}
\nonumber
\fun{K}=\fun{\esbd{B_1^+(0)}} &\le& \liminf_n \fun{\esbd{E_n}} \le \\
\nonumber
&\le&
\liminf_n \fun{H_n} \le \\
\nonumber
&\le&
\liminf_n \fun{K_n}+ c_2 \limsup_n \hs^1(L_n) \le \\
\nonumber
&\le&
\liminf_n \fun{K_n}+ 3c_2\varepsilon.
\end{eqnarray}
Letting $\varepsilon \to 0$, we have
\begin{equation}
\label{eqsemfs}
\fun{K} \le \liminf_n \fun{K_n}.
\end{equation}
\par
To obtain the thesis, we need to prove that for $\hs^1$-almost
all points $x_0$ of $K \cap U$
\begin{equation}
\label{eq3}
\limsup_{\rho \to 0^+} \frac{\mu(B_\rho(x_0))}{2\rho} \ge \varphi(x_0,\nu_{x_0})
\end{equation}
where $\nu_{x_0}$ indicates the normal to $K$ at $x_0$: this is sufficient in order to compare
$\mu_0$ and $\mu$ (see Theorem 2.56 of \cite{AFP}).
\par
Up to a rotation we may assume that $\nu_{x_0}=e_2$.
Let $\rho_k \to 0^+$ and let $T_k$ be the map defined by
$$
T_k(\xi)= \frac{1}{\rho_k} (\xi-x_0)
$$
which brings the ball $B_{\rho_k}(x_0)$ to the unit ball of the plane. By our choice of $x_0$,
$\hs^1 \res T_k(K)$ converges locally weakly star in the sense of measures to $\hs^1 \res H$ where $H$
denotes the horizontal axis of the plane.
\par
Note that for $k \to \infty$
\begin{equation}
\label{eq2}
T_k(K) \cap \Bb \to H \cap \Bb
\end{equation}
in the Hausdorff metric. In fact, up to a subsequence,
$T_k(K) \cap \Bb \to \widetilde{K}$ by compactness of the Hausdorff metric.
Clearly
$H \cap \Bb \subseteq \widetilde{K}$ because if $y \in (H \cap \Bb) \setminus
\widetilde{K}$, there exists $\rho>0$ such that $T_k(K) \cap B_{\rho}(y) =
\emptyset$ definitively and so
$$
\hs^1(H \cap \Bb \cap B_{\rho}(y)) \le \liminf_k \hs^1(T_k(K) \cap B_{\rho}(y))=0
$$
which is absurd. Conversely, $\widetilde{K} \subseteq H \cap \Bb$ because if
$y \in \widetilde{K} \setminus (H \cap \Bb)$, there exists $\rho>0$ such that $H \cap \Bb
\cap \overline{B}_{\rho}(y)=\emptyset$ and by the inequality
$$
\limsup_k \hs^1(\overline{B}_{\rho}(y) \cap T_k(K)) \le
\hs^1(\overline{B}_{\rho}(y) \cap H \cap \Bb)
$$
we deduce
\begin{equation}
\label{eqsem2}
\limsup_k \hs^1(\overline{B}_{\rho}(y) \cap T_k(K))=0.
\end{equation}
But we proved that $H \cap \Bb \subseteq \widetilde{K}$ and so the points of $H \cap \Bb$
are limit of points of $T_k(K)$: since every $T_k(K)$ is arcwise connected (they are connected
and have finite $\hs^1$ measure), we have
that $\hs^1(\overline{B}_{\rho}(y) \cap T_k(K)) \ge \rho$ definitively and this
contradicts (\ref{eqsem2}).
\par
We may suppose that $\rho_k$'s are chosen in such a way that
\begin{equation}
\label{rok}
\mu(\partial B_{\rho_k}(x_0))=0 \quad \quad
\lim_n \mu_n (\overline{B}_{\rho_k}(x_0))=\mu(B_{\rho_k}(x_0)).
\end{equation}
Since $T_k(K_{n}) \to T_k(K)$ in the Hausdorff metric for $n \to +\infty$  by (\ref{eq2}) and (\ref{rok}) there exists
a subsequence $n_k$ such that
\begin{equation}
\label{hausconv}
T_k(K_{n_k}) \cap \Bb \to H \cap \Bb
\end{equation}
in the Hausdorff metric for $k \to +\infty$ and
$$
\mu_{n_k}(\overline{B}_{\rho_k}(x_0)) \le
\mu (B_{\rho_k}(x_0))+ \rho_k^2.
$$
We now want to use the device of the first part of the proof: we employ the notation
introduced before. Let $\varepsilon>0$, $\eta>0$,
$R^\varepsilon_\eta:=S_\eta \cap V^\varepsilon$,
$\partial^{\pm}R^\varepsilon_\eta:= R^\varepsilon_\eta \cap \partial^{\pm} V^\varepsilon$;
for $k$ large enough $T_k(K_{n_k}) \cap V^\varepsilon \subseteq R^\varepsilon_\eta$ and
if $C^{\pm}_k$ is the connected component of
$(T_k(K_{n_k}) \cap \Bb) \cup \partial^{-}R^\varepsilon_\eta  \cup \partial^{+}R^\varepsilon_\eta$ containing
$\partial^{\pm}R^\varepsilon_\eta$,
we have
$(T_k(K_{n_k}) \cap \Bb) \cup \partial^{-}R^\varepsilon_\eta  \cup \partial^{+}R^\varepsilon_\eta=C^-_k \cup C^+_k$.
In fact if $\xi \notin C^-_k \cup C^+_k$ and $C_k^\xi$ be the connected component of
$(T_k(K_{n_k}) \cap \Bb) \cup \partial^{-}R^\varepsilon_\eta  \cup \partial^{+}R^\varepsilon_\eta$
containing $\xi$, by (\ref{hausconv}), $C_k^\xi \cap \partial R^\varepsilon_\eta=\emptyset$ for
$k$ large enough and
so $C_k^\xi$ would be connected against the connectedness of
$T_k(K_{n_k}) \cup \partial^{-}R^\varepsilon_\eta  \cup \partial^{+}R^\varepsilon_\eta$.
\par
By (\ref{hausconv}), we deduce easily that it is possible to join a point
of $C^+_k$ and a point of $C^-_k$ through a line $l_k \subseteq \Bb$ such that $\hs^1(l_k) \le \varepsilon$
for $k$ large enough.
\par
Considering $H_k:=(T_k(K_{n_k}) \cap \Bb) \cup L_k \cup l_k$, $H_k$ is connected in $\Bb$ and converges to
$H \cap \Bb$ in the Hausdorff metric. Applying (\ref{eqsemfs}) with $\varphi=\varphi(x_0,\cdot)$, and
since
$\sup\{|\varphi(x_0+\rho_k (\cdot), \nu) -  \varphi(x_0, \nu)|\} \to 0$ in $\Bb \times S^1$ uniformly
by the continuity of $\varphi$, we get
$$
2\varphi(x_0, e_2) \le
\liminf_k \int_{H_k} \varphi(x_0, \nu) \,d\hs^1 \le
\liminf_k \int_{T_k(K_{n_k}) \cap \Bb} \varphi(x_0+\rho_k x,\nu) \,d\hs^1
+3c_2\varepsilon.
$$
Letting $\varepsilon \to 0$, we obtain
\begin{equation}
\label{eqsemfs2}
2\varphi(x_0, e_2) \le
\liminf_k \int_{T_k(K_{n_k}) \cap \Bb} \varphi(x_0+\rho_k x,\nu) \,d\hs^1.
\end{equation}
\par
Now we are ready to conclude: in fact
\begin{eqnarray}
\nonumber
\limsup_{\rho \to 0} \frac{\mu(B_{\rho}(x_0))}{2\rho} &\ge&
\limsup_k \frac{\mu(B_{\rho_k}(x_0))}{2\rho_k} \ge \\
\nonumber
&\ge&
\liminf_k \frac{\mu_{n_k}(\overline{B}_{\rho_k}(x_0))}{2\rho_k} =\\
\nonumber
&=&
\frac{1}{2} \liminf_k \int_{T_k(K_{n_k}) \cap \Bb} \varphi(x_0+\rho_k x, \nu) \,d\hs^1
\ge \varphi(x_0,e_2),
\end{eqnarray}
the last inequality coming from (\ref{eqsemfs2}).
\par
Let's now turn to the case $m \ge 2$.
Let $(K_n) \in \ks_m^f(\Omb)$ converges to $K$; up to a
subsequence, we may suppose that there exists $m' \le m$ such that each $K_n$
has exactly $m'$ connected components $\widehat{K}_n^1, \cdots,\widehat{K}_n^{m'}$.
We may suppose moreover that for all $i$, $\widehat{K}_n^{i} \to \widehat{K}^{i}$ in the Hausdorff
metric: it is readily seen that $K = \cup_{i=1}^{m'} \widehat{K}^{i}$ so that, using
the lower semicontinuity for the case $m=1$, we obtain
\begin{eqnarray}
\nonumber
\funU{K} &\le& \sum_{i=1}^{m'} \funU{\widehat{K}^{i}} \le \\
\nonumber
&\le& \liminf_n \sum_{i=1}^{m'} \funU{\widehat{K}_n^{i}} = \\
\nonumber
&=& \liminf_n \funU{K_n}.
\end{eqnarray}
\end{proof}

Theorem \ref{main0} is now proved: it is sufficient to apply Theorem \ref{sem.2} with
$U=B_R(0)$, $\Omb \subseteq B_R(0)$.

\begin{corollary}
\label{sem.4}
Let $(H_n)$ be a sequence in $\ks(\Omb)$ which converges to $H$ in the
Hausdorff metric. Let $m \ge 1$ and let $(K_n)$ be a sequence in $\ks_m^f(\Omb)$
which converges to $K$ in the Hausdorff metric. Then
$$
\fun{K \setminus H} \le \liminf_n \fun{K_n \setminus H_n}.
$$
\end{corollary}

\begin{proof}
Let $\varepsilon>0$ and let $H^\varepsilon= \{x \in \Omb\,:\, {\rm
dist}(x,H)\le \varepsilon\}$. Definitively $H_n \subseteq H^\varepsilon$ so that
$K_n \setminus H^\varepsilon \subseteq K_n \setminus H_n$. Applying Theorem
\ref{sem.2} with $U= \R^2 \setminus H^\varepsilon$, we have
$$
\fun{K \setminus H^\varepsilon} \le \liminf_n \fun{K_n \setminus H^\varepsilon}
\le \liminf_n \fun{K_n \setminus H_n}.
$$
Letting $\varepsilon$ go to zero, we obtain the thesis.
\end{proof}

The following result will be useful in sections \ref{antplan} and \ref{plan}.

\begin{theorem}
\label{sem.5}
Given $m \ge 1$, let $(H_n)$ be a sequence in $\ks_m^f(\Omb)$ which converges to $H$ in the Hausdorff metric,
and let $K \in \ks_m^f(\Omb)$ with $H \subseteq K$. Then there exists a sequence $(K_n)$ in
$\ks_m^f(\Omb)$ which converges to $K$ in the Hausdorff metric and such that
$H_n \subseteq K_n$ and
\begin{equation}
\label{eqs1}
\fun{K \setminus H} = \lim_n \fun{K_n \setminus H_n}.
\end{equation}
\end{theorem}

\begin{proof}
Following Lemma 3.8 of \cite{DMT}, the connected components $C_i$ of $K \setminus H$ are at least countable and
satisfy $\hs^1(\overline{C_i})=\hs^1(C_i)$.
Since $H_n \to H$ in the Hausdorff metric and $\Om$ has Lipschitz boundary, we can find arcs
$Z_n^i$ in $\Omb$ joining $H_n$ and $C_i$ such that $\hs^1(Z_n^i) \to 0$ as $n \to \infty$.
Given $h$, consider $K_n^h:=\cup_{i=1}^{h} Z_n^i \cup \cup_{i=1}^{h}\overline{C}_i$; we have
$K_n^h \in \ks^f_h(\Omb)$, $K_n^h \to K^h:=\cup_{i=1}^{h}\overline{C}_i$ in the Hausdorff metric. Note that
$\nu \hs^1 \res K_n^h \to \nu \hs^1 \res K^h$ strictly for $n \to \infty$. By Theorem \ref{Resh2}, since
$\hs^1(\overline{C_i})=\hs^1(C_i)$,
we have
\begin{eqnarray}
\nonumber
\lim_n \fun{K_n^h} &=& \fun{K^h} = \\
\nonumber
&=& \fun{\cup_{i=1}^h C_i} \le \fun{K \setminus H}.
\end{eqnarray}
Choose $h_n \to +\infty$ such that
$$
\limsup_n \fun{K_n^{h_n}} \le \fun{K \setminus H},
$$
so that
$$
\lim_n \sum_{i=1}^{h_n} \hs^1(Z_n^i)=0.
$$
If we pose $K_n:= H_n \cup K_n^{h_n}$, we have
$K_n \in \ks^f_m(\Omb)$, $K_n \to K$ in the Hausdorff metric and
$$
\limsup_n \fun{K_n \setminus H_n} \le
\limsup_n \fun{K_n^{h_n}} \le \fun{K \setminus H}.
$$
The converse inequality comes from Corollary \ref{sem.4}.
\end{proof}

\section{THE ANTI-PLANAR ANISOTROPIC CASE}
\label{antplan}

In this section we deal with quasi-static growth of brittle fractures in
inhomogeneous anisotropic linearly elastic bodies under
anti-planar displacements. We employ the notation of Section \ref{main}.
\par
We begin with the following lemma which extends Theorem \ref{p2} considering boundary data.
The idea is due to A. Chambolle.

\begin{lemma}
\label{approx1}
Let $m \ge 1$, $K_n$ a sequence in $\ks_m(\Omb)$ which
converges to $K$ in the Hausdorff metric and such that
$\leb^2(\Om \setminus K_n) \to \leb^2(\Om \setminus K)$. Let $g_n \to g$
strongly in $H^1(\Om)$ and let
$\Gamma(g_n,K_n)$ and $\Gamma(g,K)$ be the sets introduced in (\ref{gamma}).
Then for every $u \in \Gamma(g,K)$,
there exists $u_n \in \Gamma(g_n,K_n)$ such that $\nabla u_n \to \nabla u$ strongly
in $L^2(\Om,\R^2)$.
\end{lemma}

\begin{proof}
Consider $\Om'$ a regular open set containing $\overline{\Om}$ and pose
$\partial_N \Om := \partial \Om \setminus \partial_D \Om$. Since $\partial \Om$ is regular,
we may extend $g_n$ and  $g$ to $H^1(\Om')$ and suppose $g_n \to g$ strongly in $H^1(\Om')$.
Note that if $H_n:=K_n \cup \overline{\partial_N \Om}$ and $H= K \cup \overline{\partial_N \Om}$,
$H_n,H \in \ks_{m'}(\overline{\Om'})$, $H_n \to H$ in the Hausdorff metric and
$\leb^2(\Om' \setminus H_n) \to \leb^2(\Om' \setminus H)$. Consider
$$
v:= \left\{
\begin{array}{l}
u \\ g
\end{array}
\begin{array}{l}
{\rm in}\; \Om \\
{\rm in}\; \Om' \setminus \Om
\end{array}
\right.
$$
Clearly $v \in L^{1,2}(\Om' \setminus H)$; we may apply Theorem \ref{p2} and deduce that
there exists $v_n \in L^{1,2}(\Om' \setminus H_n)$ such that $\nabla v_n \to \nabla v$
strongly in $L^2(\Om',\R^2)$. Note that we may assume $(v_n-v)$ has null average on $\Om' \setminus \Omb$,
because we are allowed to add constants to $v_n$; since $\Om' \setminus \Omb$ is regular, by Poincar\'e
inequality we obtain $v_n \to v$ strongly in $H^1(\Om' \setminus \Omb)$.
Let $E_\Om$ be a linear extension operator from $H^1(\Om' \setminus \Omb)$ to $H^1(\Om')$. If
$w_n:= (v_n-v)|_{\Om' \setminus \Omb}$, we can choose
$$
u_n:= v_n -E_\Om w_n +(g_n-g)
$$
restricted to $\Om$. It is readily seen that $u_n \in \Gamma(g_n,K_n)$ and $\nabla u_n \to \nabla u$
strongly in $L^2(\Om,\R^2)$.
\end{proof}

By standard arguments, it can be proved that the minimum of problem (\ref{prantplan}) is attained.
Moreover, it can be
shown that, since $\Om \setminus K$ is not guaranteed to be regular, this minimum is in general not attained
in $H^1(\Om \setminus K)$ when the boundary data $g$ is not bounded: the reader is referred to \cite{Ma}.
The following proposition deals with the behavior of minima when the compact set $K$ varies.

\begin{proposition}
\label{antplan0}
Let $m \ge 1$, $\lambda \ge 0$, $(K_n)$ a sequence in $\ks_m^\lambda(\Omb)$ which
converges to $K$ in the Hausdorff metric, $(g_n)$ a sequence in $H^1(\Om)$ which
converges to $g$ strongly in $H^1(\Om)$. Let $u_n$ be a solution of the minimum problem
\begin{equation}
\label{eqantplan01}
\min_{v \in \Gamma(g_n,K_n)} \nL{u}^2
\end{equation}
and let $u$ be a solution of the minimum problem
\begin{equation}
\label{eqantplan02}
\min_{v \in \Gamma(g,K)} \nL{u}^2,
\end{equation}
where $\Gamma(g_n,K_n)$ and $\Gamma(g,K)$ are defined as in (\ref{gamma}).
\par
Then $\nabla u_n \to \nabla u$ strongly in $L^2(\Om, \R^2)$.
\end{proposition}

\begin{proof}
Using $g_n$ as test function, we obtain
$$
\nL{u_n} \le \nL{g_n} \le c <+\infty.
$$
By (\ref{unifellipt}), there exists $\psi \in L^2(\Om,\R^2)$ such that, up to a subsequence,
$\nabla u_n \weak \psi$ weakly in $L^2(\Om,\R^2)$.
It is not difficult to prove that there exists $u \in L^2_{loc}(\Om)$ such that $\nabla u=\psi$ in $\Om \setminus K$.
Moreover by means of Poincar\'e inequality, we deduce that $u=g$ on $\partial_D \Om \setminus K$.
According to Lemma \ref{approx1}, let $v_n \in \Gamma(g_n,K_n)$ with $\nabla v_n \to \nabla u$ strongly in $L^2(\Om,\R^2)$;
since $||u_n||_a \le ||v_n||_a$ by minimality of $u_n$, we obtain
$\limsup_n ||u_n||_a \le ||u||_a$. This proves $\nabla u_n \to \nabla u$ strongly in $L^2(\Om, \R^2)$.
\end{proof}

We now turn to the proof of Theorem \ref{main1}. We use a discretization in time.
Given $\delta>0$, let $N_\delta$ be the largest integer such that
$\delta N_\delta \le 1$; for $i \ge 0$ we pose $t_i^\delta=i\delta$ and for
$0 \le i \le N_\delta$ we pose $g_i^\delta=g(t_i^\delta)$. Define $K_i^\delta$ as a solution
of the minimum problem
\begin{equation}
\label{eqac2}
\min_K \left\{ \eub(g_i^\delta,K) \,:\,
K \in \ks_m^f(\Omb),\, K_{i-1}^\delta \subseteq K \right\},
\end{equation}
where $K_{-1}^\delta=K_0$.

\begin{lemma}
\label{ac.2}
The minimum problem (\ref{eqac2}) admits a solution.
\end{lemma}

\begin{proof}
We proceed by induction. Suppose $K_{i-1}^\delta$ is constructed and that
$\lambda > \eub(g_i^\delta,K_{i-1}^\delta)$. Let $(K_n)$ be a minimizing sequence of problem
(\ref{eqac2}) and let $u_n$ be a solution of the minimum problem (\ref{prantplan}) which defines
$\eub(g_i^\delta,K_n)$. Up to a subsequence, $K_n \to K$ in the Hausdorff metric and
$K_{i-1}^\delta \subseteq K$.
Since
$$
\nL{u_n}^2 + \fun{K_n} \le \lambda
$$
for $n$ large, we have that
$$
\fun{K_n} \le \lambda.
$$
We have  $K_n \in \ks^{\alpha_1^{-1} \lambda}_m(\Omb)$ and applying Proposition \ref{antplan0}, we have
$||u_n||_a \to ||u||_a$
where $u$ is a solution of problem (\ref{prantplan}) which defines $\eub(g_i^\delta,K)$;
moreover by Theorem \ref{sem.2}, we get
$$
\fun{K} \le \liminf_n \fun{K_n} \le \lambda.
$$
Thus $K \in \ks^f_m(\Omb)$ and $\eub(g_i^\delta,K) \le \liminf_n \eub(g_i^\delta,K_n)$.
We conclude that $K$ is a solution of the minimum problem (\ref{eqac2}).
\end{proof}

Now, consider the following piecewise constant interpolation:
put $g^\delta(t)=g_i^\delta$, $K^\delta(t)=K_i^\delta$, $u^\delta(t)=u_i^\delta$ for
$t_i^\delta \le t < t_{i+1}^\delta$, where $u_i^\delta$ is a solution of problem
(\ref{prantplan}) which defines $\eub(g_i^\delta,K_i^\delta)$.

\begin{lemma}
\label{ac.3}
There exists a positive function $\rho(\delta)$, converging to zero as $\delta \to 0$, such
that for all $s<t$ in $[0,1]$,
\begin{eqnarray}
\label{eqac3}
\nL{u^\delta(t)}^2 + \fun{K^\delta(t)} &\le&
\nL{u^\delta(s)}^2 + \fun{K^\delta(s)} +\\
\nonumber
&&+2 \int_{t_i^\delta}^{t_j^\delta} \pr{u^\delta (t)}{\dot{g}(t)}\,dt +\rho(\delta)
\end{eqnarray}
where $t_i^\delta \le s < t_{i+1}^\delta$ and $t_j^\delta \le t < t_{j+1}^\delta$.
\end{lemma}

\begin{proof}
Inequality (\ref{eqac3}) is precisely
$$
\nL{u_j^\delta}^2 + \fun{K_j^\delta} \le
\nL{u_i^\delta}^2 + \fun{K_i^\delta}
+2 \int_{t_i^\delta}^{t_j^\delta} \pr{u^\delta (t)}{\dot{g}(t)}\,dt +\rho(\delta).
$$
To obtain this one, it is sufficient to adapt the proof of Lemma 7.3 in \cite{DMT}.
\end{proof}

\begin{lemma}
\label{ac4}
There exists a constant $C$, depending only on $g$ and $K_0$, such that
$$
\nL{u^\delta(t)} \le C \quad \quad \fun{K^\delta(t)} \le C
$$
for every $\delta>0$ and $t \in [0,1]$.
In particular, there exists $\lambda >0$ such that for all $t \in [0,1]$,
$K^\delta(t) \in \ks^\lambda_m(\Omb)$.
\end{lemma}

\begin{proof}
Put $\eta=\max_t\{\nL{g(t)}, \nL{\dot{g}(t)}\}$. Clearly
$\nL{u^\delta(t)} \le \nL{g^\delta(t)}
\le \eta$ since $g^\delta(t)$
is an admissible displacement for $K^\delta(t)$.
Clearly from inequality (\ref{eqac3}) with $s=0$, we obtain
\begin{eqnarray}
\nonumber
\nL{u^\delta(t)}^2 + \fun{K^\delta(t)} &\le&
\nL{u^\delta(0)}^2 + \fun{K^\delta(0)} +\\
\nonumber
&& +2 \int_{0}^{t_j^\delta} \pr{u^\delta (t)}{\dot{g}(t)}\,dt +\rho(\delta) \le  \\
\nonumber
&\le& \nL{u_0^\delta}^2 + \fun{K_0}+
2\eta^2 +\rho(\delta).
\end{eqnarray}
The last term depends only on $g$ and $K_0$ and so we obtain the first part of the thesis.
The second one comes from (\ref{varphiprop}).
\end{proof}

\begin{lemma}
\label{ac5}
Let $C$ be the constant of Lemma \ref{ac4}. There exists an increasing function
$K:[0,1] \to \ks_m^f(\Omb)$ (that is $K(s) \subseteq K(t)$ for every $0 \le s \le t \le 1$),
such that, for every $t \in [0,1]$, $K^\delta(t)$ converges to $K(t)$ in the Hausdorff
metric as $\delta \to 0$ along a suitable sequence independent of $t$. Moreover if $u(t)$ is
a solution of the minimum problem (\ref{prantplan}) which defines $\eub(g(t),K(t))$, for
every $t \in [0,1]$ we have $\nabla u^\delta(t) \to \nabla u(t)$ strongly in
$L^2(\Om,\R^2)$.
\end{lemma}

\begin{proof}
The first part is a variant of Helly's theorem for monotone function: for a proof see Lemma 7.5 of
\cite{DMT}; the second part comes directly from Lemma \ref{ac4} and Proposition \ref{antplan0}.
\end{proof}

Fix now the sequence $(\delta_n)$ and the increasing map $t \to K(t)$ given by Lemma \ref{ac5}.
We indicate $K^{\delta_n}(t)$ by $K_n(t)$ and $u^{\delta_n}(t)$ by $u_n(t)$.
\par
The following property of the pair $(g(t), K(t))$ is important for subsequent results.

\begin{lemma}
\label{ac6}
For every $t \in [0,1]$ we have
\begin{equation}
\label{eqac4}
\eub(g(t),K(t)) \le \eub(g(t),K) \quad \forall K \in \ks_m^f(\Omb),\, K(t) \subseteq K.
\end{equation}
Moreover
\begin{equation}
\label{eqac5}
\eub(g(0),K(0)) \le \eub(g(0),K) \quad \forall K \in \ks_m^f(\Omb),\, K_0 \subseteq K.
\end{equation}
\end{lemma}

\begin{proof}
Let $t \in [0,1]$ and $K \in \ks_m^f(\Omb)$ with $K(t) \subseteq K$. Since $K_n(t) \to
K(t)$ in the Hausdorff metric as $\delta_n \to 0$, by Theorem \ref{sem.5} there exists a
sequence $(K_n)$ in $\ks_m^f(\Omb)$ converging to $K$ in the Hausdorff metric, such
that $K_n(t) \subseteq K_n$ and
\begin{equation}
\label{eqac6}
\fun{K_n \setminus K_n(t)} \to
\fun{K \setminus K(t)}.
\end{equation}
By Lemma \ref{ac4}, there exists $\lambda>0$ such that $K_n(t) \in \ks^\lambda_m(\Omb)$ for
all $n$. By (\ref{eqac6}), we deduce that there exists $\lambda'>\lambda$
with $K_n \in \ks^{\lambda'}_m(\Omb)$ for all $n$.
\par
Let $v_n$ and $v$ solutions of problems (\ref{prantplan}) which define
$\eub(g_n(t),K_n)$ and $\eub(g(t),K)$. By minimality of $K_n(t)$ we have
$\eub(g_n(t),K_n(t)) \le \eub(g_n(t),K_n)$ and so
\begin{equation}
\label{ineq}
\nL{u_n(t)}^2 \le \nL{v_n}^2 + \fun{K_n \setminus K_n(t)};
\end{equation}
as $\delta_n \to 0$, $\nabla u_n(t) \to \nabla u(t)$ and $\nabla v_n \to \nabla v$ strongly in
$L^2(\Om,\R^2)$ by Proposition \ref{antplan0}: passing to the limit in (\ref{ineq}) and
adding to both sides $\fun{K(t)}$, by (\ref{eqac6}) we have the thesis.
\par
A similar proof holds for (\ref{eqac5}).
\end{proof}

\begin{lemma}
\label{ac9}
The function $t \to \eub(g(t),K(t))$ is absolutely continuous and
$$
\frac{d}{dt} \eub(g(t),K(t))=2\pr{u(t)}{\dot{g}(t)}
\quad {\it for\,a.e}\; t\in [0,1]
$$
where $u(t)$ is a solution of the minimum problem (\ref{prantplan}) which defines
$\eub(g(t),K(t))$.
\end{lemma}

\begin{proof}
We rewrite (\ref{eqac3}) in the following form
$$
\nL{u_n(t)}^2 + \fun{K_n(t) \setminus K_n(s)} \le
\nL{u_n(s)}^2
+2 \int_{t_i^{\delta_n}}^{t_j^{\delta_n}} \pr{u_n (t)}{\dot{g}(t)}\,dt +\rho(\delta_n)
$$
for $s\le t$ and $t_i^{\delta_n} \le s < t_{i+1}^{\delta_n}$ and $t_j^{\delta_n} \le t < t_{j+1}^{\delta_n}$.
Passing to the limit for $\delta_n \to 0$, using Corollary \ref{sem.4} we obtain
$$
\nL{u(t)}^2 + \fun{K(t) \setminus K(s)} \le
\nL{u(s)}^2
+2 \int_{s}^{t} \pr{u(\tau)}{\dot{g}(\tau)}\,d\tau,
$$
so that
\begin{eqnarray}
\nonumber
\nL{u(t)}^2 + \fun{K(t)} &\le&
\nL{u(s)}^2 + \fun{K(s)} +\\
\nonumber
&& +2 \int_{s}^{t} \pr{u(\tau)}{\dot{g}(\tau)}\,d\tau.
\end{eqnarray}
Following Lemma 6.5 of \cite{DMT}, we can prove that the function $F(g):= \eub(g,K(t))$ is
differentiable on $H^1(\Om)$ and its differential is given by $dF(g)h= 2\pr{u(t)}{h}$ where $u(t)$
is a solution of problem (\ref{prantplan}) which defines $\eub(g,K(t))$.
By Lemma \ref{ac6}, we obtain
$$
\eub(g(t),K(t)) -\eub(g(s),K(s)) \ge
\eub(g(t),K(t)) -\eub(g(s),K(t)) =
2\int_s^t \pr{u(\tau,t)}{\dot{g}(\tau)} d\tau
$$
where $u(\tau,t)$ is a solution of the minimum problem (\ref{prantplan})
which defines $\eub(g(\tau),K(t))$.
We can conclude that $t \to \eub(g(t),K(t))$ is absolutely continuous since $\nL{u(t)}$ and
$\nL{u(\tau,t)}$ are bounded by Lemma \ref{ac4}. Moreover, dividing the previous inequalities
by $t-s$ and letting $s \to t$, since $\nabla u(\tau,t) \to \nabla u(t)$ strongly in $L^2(\Om,\R^2)$ for
$\tau \to t$, we obtain
$$
\frac{d}{dt} \eub(g(t),K(t))=2\pr{u(t)}{\dot{g}(t)}
\quad {\it for\,a.e}\; t\in [0,1].
$$
\end{proof}

We now turn to the proof of Theorem \ref{main1}.
Points $(a)$ and $(b)$ are proved in lemmas \ref{ac5} and \ref{ac6} while points $(d)$ and $(e)$ are proved in Lemma \ref{ac9}.
Point $(f)$ and its equivalence to point $(e)$ stated in Remark \ref{mainrem} are proved adapting Lemma 6.4 of \cite{DMT}.
To prove point $(c)$, we need the following lemma.

\begin{lemma}
\label{ac8}
Let $K\,:\,[0,1] \to \ks_m^f(\Omb)$ be a map which satisfies lemmas \ref{ac6} and \ref{ac9}. Then for every
$t \in ]0,1]$,
$$
\eub(g(t),K(t)) \le \eub(g(t),K) \quad \forall K \in \ks_m^f(\Omb)\,: \cup_{s<t} K(s) \subseteq K.
$$
\end{lemma}

\begin{proof}
Consider $t \in ]0,1]$ and $K \in \ks_m^f(\Omb)$ such that $\cup_{s<t} K(s) \subseteq K$. For $0\le s < t$ we have $K(s) \subseteq K$
and so by Lemma \ref{ac6}, $\eub(g(s),K(s)) \le \eub(g(s),K)$. By Lemma \ref{ac9}, these expressions continuously depend
on $s$ and so passing to the limit for $s \to t$, we obtain the thesis.
\end{proof}

Consider now the particular case in which $g(0)=0$: there exists a solution $K(t)$ to the problem of
evolution such that $K(0)=K_0$ because in the time discretization method employed, we can choose $K^\delta(0)=K_0$.
Under this assumption, we prove that this method gives an approximation of the energy of the solution.
\par
We pose
$$
\eub_n(t)=\nL{u_n(t)}^2 +\fun{K_n(t)}
$$
and
$$
\eub(t)=\eub(g(t),K(t))=\nL{u(t)}^2 +\fun{K(t)}.
$$
The following convergence result holds.

\begin{proposition}
\label{ac7}
For all $t \in [0,1]$ the following facts hold:
\begin{enumerate}
\item[(a)]
$K_n(t) \to K(t)$ in the Hausdorff metric;
\item[(b)]
$\nabla u_n(t) \to \nabla u(t)$ strongly in $L^2(\Om,\R^2)$;
\item[(c)]
$\fun{K_n(t)} \to \fun{K(t)}.$
\end{enumerate}
In particular $\eub_n(t) \to \eub(t)$ for all $t \in [0,1]$.
\end{proposition}

\begin{proof}
We have already proved points $(a)$ and $(b)$ in Lemma \ref{ac5}.
Since the functions $t \to \fun{K_n(t)}$ are increasing and bounded, we may suppose that,
by Helly's theorem, they converge pointwise to a bounded
increasing function $h:\,[0,1] \to [0,\infty[$ i.e. for all $t \in [0,1]$
$$
\lim_n \fun{K_n(t)}= h(t).
$$
Moreover by Theorem \ref{sem.2} we have that $\fun{K(t)} \le h(t)$ for all $t \in [0,1]$
and by construction $\lambda(0)=\fun{K(0)}$;
in particular we have for all $t \in [0,1]$
$$
\eub(t) \le \nL{u(t)}^2 +h(t)
$$
and $\eub(0)=\nL{u(0)}^2+h(0)$.
Passing to the limit in (\ref{eqac3}), by $(b)$ we obtain
\begin{equation}
\label{aceqlam}
\nL{u(t)}^2 + h(t) \le
\nL{u(s)}^2 + h(s)
+2 \int_{s}^{t} \pr{u (t)}{\dot{g}(t)}\,dt.
\end{equation}
Since by condition ${\rm (e)}$ of Theorem \ref{main1}
$$
\eub(t)-\eub(0)= 2 \int_0^t \pr{u(\tau)}{\dot{g}(\tau)}\,d\tau,
$$
we have
\begin{eqnarray}
\nonumber
\nL{u(t)}^2 +h(t) &-& \eub(t)= \\
\nonumber
&\le&
2\int_0^t \pr{u(\tau)}{\dot{g}(\tau)}\,d\tau - 2\int_0^t \pr{u(\tau)}{\dot{g}(\tau)}\,d\tau =0.
\end{eqnarray}
We conclude that $h(t)=\fun{K(t)}$ for all $t \in [0,1]$. This proves point $(c)$ and
the thesis is obtained.
\end{proof}

\section{THE PLANAR ANISOTROPIC CASE}
\label{plan}

In this section we breafily sketch the modifications of the arguments used in the previous section
in order to deal with the evolution of fractures in
inhomogeneous anisotropic
linearly elastic bodies under planar displacements.
We employ the notation of Section \ref{main}.
\par
The following lemma can be obtained
with arguments similar to those of Lemma \ref{approx1}.

\begin{lemma}
\label{approx2}
Let $m \ge 1$, $K_n$ a sequence in $\ks_m(\Omb)$ which
converges to $K$ in the Hausdorff metric and such that
$\leb^2(\Om \setminus K_n) \to \leb^2(\Om \setminus K)$. Let $g_n \to g$
strongly in $H^1(\Om, \R^2)$ and let
$\vub(g_n,K_n)$ and $\vub(g,K)$ be the sets introduced in (\ref{ni}).
Then for every $u \in \vub(g,K)$,
there exists $u_n \in \vub(g_n,K_n)$ such that $Eu_n \to Eu$ strongly
in $L^2(\Om,\msim)$.
\end{lemma}

By standard techniques, it can be proved that the minimum in problem
(\ref{prplan}) is attained. The following result is similar to Proposition \ref{antplan0}
and deals with the behavior of these minima as $K$ varies.

\begin{proposition}
\label{plan0}
Let $m \ge 1$ and $\lambda \ge 0$, let $K_n$ be a sequence in $\ks_m^\lambda(\Omb)$ which
converges to $K$ in the Hausdorff metric, and let $g_n$ be a sequence in $H^1(\Om)$ which
converges to $g$ strongly in $H^1(\Om)$. Let $u_n$ be a solution of the minimum problem
\begin{equation}
\label{eqplan01}
\min_{v \in \vub(g_n,K_n)} \nA{v}^2,
\end{equation}
and let $u$ be a solution of the minimum problem
\begin{equation}
\label{eqplan02}
\min_{v \in \vub(g,K)} \nA{v}^2
\end{equation}
where $\vub(g_n,K_n)$ and $\vub(g,K)$ are defined as in (\ref{ni}).
\par
Then $Eu_n \to Eu$ strongly in $L^2(\Om,\msim)$.
\end{proposition}

\begin{proof}
Using $g_n$ as test function we obtain $\nA{u_n} \le \nA{g_n} \le c <+\infty$.
By assumption on $A$, there exists $\sigma \in L^2(\Om,\msim)$ such that up to a subsequence $Eu_n \weak \sigma$ weakly in
$L^2(\Om,\msim)$.
It is not difficult to prove that there exists $u \in L^2_{loc}(\Om,\R^2)$ such that $Eu=\sigma$ in $\Om \setminus K$.
Moreover by means of Korn-Poincar\'e inequality, we deduce that $u=g$ on $\partial_D \Om \setminus K$.
According to Lemma \ref{approx2}, let $v_n \in \vub(g_n,K_n)$ with $Ev_n \to Eu$ strongly in $L^2(\Om,\msim)$;
since $\nA{u_n} \le \nA{v_n}$ by minimality of $u_n$, we obtain
$$
\limsup_n \nA{u_n} \le \limsup_n \nA{v_n}= \nA{u}.
$$
This proves $Eu_n \to Eu$ strongly in $L^2(\Om,\msim)$.
\end{proof}

We employ again a time discretization process. As before
given $\delta>0$, let $N_\delta$ be the largest integer such that
$\delta N_\delta \le 1$; for $i \ge 0$ we pose $t_i^\delta=i\delta$ and for
$0 \le i \le N_\delta$ we pose $g_i^\delta=g(t_i^\delta)$. Define $K_i^\delta$ as a solution
of the minimum problem
\begin{equation}
\label{eqplan1}
\min_K \left\{ \gs(g_i^\delta,K) \,:\,
K \in \ks_m^f(\Omb),\, K_{i-1}^\delta \subseteq K \right\},
\end{equation}
where $K_{-1}^\delta=K_0$.

\begin{lemma}
\label{plan1}
The minimum problem (\ref{eqplan1}) admits a solution.
\end{lemma}

\begin{proof}
We proceed by induction. Suppose $K_{i-1}^\delta$ is constructed and that
$\lambda > \gs(g_i^\delta,K_{i-1}^\delta)$. Let $(K_n)$ be a minimizing sequence of problem
(\ref{eqplan1}) and let $u_n$ be a solution of the minimum problem (\ref{prplan}) which defines
$\gs(g_i^\delta,K_n)$. Up to a subsequence $K_n \to K$ in the Hausdorff metric and
$K_{i-1}^\delta \subseteq K$.
Since
$$
\nA{u_n}^2+ \fun{K_n} \le \lambda
$$
for $n$ large enough, we have that
$$
\fun{K_n} \le \lambda;
$$
We have $K_n \in \ks^{\alpha_1^{-1} \lambda}_m(\Omb)$ and applying Proposition \ref{plan0}, we have
$\nA{u_n} \to \nA{u}$
where $u$ is a solution of problem (\ref{prplan}) which defines $\eub(g_i^\delta,K)$;
by Theorem \ref{sem.2}, we get
$$
\fun{K} \le \liminf_n \fun{K_n} \le \lambda.
$$
Thus $K \in \ks^f_m(\Omb)$ and $\gs(g_i^\delta,K) \le \liminf_n \gs(g_i^\delta,K_n)$.
We conclude that $K$ is a solution of the minimum problem (\ref{eqplan1}).
\end{proof}

Consider as before the piecewise constant interpolation
obtained putting $g^\delta(t)=g_i^\delta$, $K^\delta(t)=K_i^\delta$, $u^\delta(t)=u_i^\delta$ for
$t_i^\delta \le t < t_{i+1}^\delta$, where $u_i^\delta$ is a solution of problem
(\ref{prplan}) which defines $\gs(g_i^\delta,K_i^\delta)$.

\begin{lemma}
\label{plan2}
There exists a positive function $\rho(\delta)$, converging to zero as $\delta \to 0$, such
that for all $s<t$ in $[0,1]$
\begin{eqnarray}
\label{eqpl3}
\nonumber
\nA{u^\delta(t)}^2 + \fun{K^\delta(t)} &\le&
\nA{u^\delta(s)}^2 + \fun{K^\delta(s)} +\\
\nonumber
&&+2 \int_{t_i^\delta}^{t_j^\delta} \prA{u^\delta (t)}{\dot{g}(t)}\,dt +\rho(\delta)
\end{eqnarray}
where $t_i^\delta \le s < t_{i+1}^\delta$ and $t_j^\delta \le t < t_{j+1}^\delta$. In particular
there exists $C>0$ depending only on $g$ and $K_0$ such that for all $t \in [0,1]$
$$
\nA{u^\delta(t)} \le C  \quad \quad \fun{K^\delta(t)} \le C.
$$
\end{lemma}

\begin{proof}
It is sufficient to adapt lemmas \ref{ac.3} and \ref{ac4}.
\end{proof}

Using Proposition \ref{plan0} and the previous lemma we obtain

\begin{lemma}
\label{plan3}
There exists an increasing function
$K:[0,1] \to \ks_m^f(\Omb)$ (that is $K(s) \subseteq K(t)$ for every $0 \le s \le t \le 1$),
such that, for every $t \in [0,1]$, $K^\delta(t)$ converges to $K(t)$ in the Hausdorff
metric as $\delta \to 0$ along a suitable sequence independent of $t$. Moreover if $u(t)$ is
a solution of the minimum problem (\ref{prplan}) which defines $\gs(g(t),K(t))$, for
every $t \in [0,1]$ we have $Eu^\delta(t) \to Eu(t)$ strongly in $L^2(\Om,\msim)$.
\end{lemma}

The proof of Theorem \ref{main2} can now be obtained using arguments similar to those of lemmas
\ref{ac6}, \ref{ac9} and \ref{ac8} of Section \ref{antplan}.
\par
Consider now the particular case in which $g(0)=0$: there exists a solution $K(t)$ to the problem of
evolution such that $K(0)=K_0$ because in the time discretization method employed we can
choose $K^\delta(0)=K_0$. Under this assumption, as in the anti-planar case, the discretization method
gives an approximation of the energy of the solution.
\par
In fact, if we pose $K_n(t):=K^{\delta_n}(t)$ and
$$
\gs_n(t):= \nA{u_n(t)}^2 +\fun{K_n(t)},
$$
$$
\gs(t):= \gs(g(t),K(t))=\nA{u(t)}^2 +\fun{K(t)},
$$
the following approximation result holds.

\begin{proposition}
\label{plan5}
As $\delta_n \to 0$ for all $t \in [0,1]$ the following facts hold:
\begin{enumerate}
\item[(a)]
$K_n(t) \to K(t)$ in the Hausdorff metric;
\item[(b)]
$Eu_n(t) \to Eu(t)$ strongly in $L^2(\Om,\msim)$;
\item[(c)]
$\fun{K_n(t)} \to \fun{K(t)}$.
\end{enumerate}
In particular $\gs_n(t) \to \gs(t)$ for all $t \in [0,1]$.
\end{proposition}

\begin{proof}
It is sufficient to adapt Proposition \ref{ac7}.
\end{proof}

\bigskip
\bigskip
\centerline{ACKNOWLEDGEMENTS}
\bigskip
The author wishes to thank Gianni Dal Maso for having proposed him
the problem, and for many helpful and interesting discussions.

\end{document}